\documentclass[12pt]{amsart}
\usepackage{graphicx}
\usepackage{amssymb,amstext,amsfonts,amscd}

\newtheorem{theorem}{Theorem}
\newtheorem{proposition}[theorem]{Proposition}
\newtheorem{lemma}[theorem]{Lemma}

\theoremstyle{definition} %%%%%
\newtheorem{definition}[theorem]{Definition}
\newtheorem{remark}[theorem]{Remark}

\numberwithin{equation}{section}
\numberwithin{theorem}{section}

\newcommand{\bz}{\mathbb{Z}}
\newcommand{\br}{\mathbb{R}}
\newcommand{\bc}{\mathbb{C}}
\newcommand{\bp}{\mathbb{P}}

\newcommand{\bff}{\mathbb{F}}   %%% Hirzebruch surfaces
\newcommand{\bl}{\mathbb{L}}    %%% lattice
\newcommand{\bu}{\mathbb{U}}    %%% hyperbolic plane
\newcommand{\M}{\mathcal{M}}
\newcommand{\Aa}{\mathcal{A}}
\newcommand{\Ea}{\mathcal{E}}
\newcommand{\Fa}{\mathcal{F}}

\newcommand{\La}{\mathcal{L}}

\newcommand{\da}{\mathcal{D}}

\newcommand{\bv}{\mathbf{v}}

\newcommand{\zero}{\mathbf{0}}

\newcommand{\rk}{\mathop\mathrm{rank}}
\newcommand{\Pic}{\mathop\mathrm{Pic}}
\newcommand{\BL}{\mathrm{bl}}
\newcommand{\mmod}{\ \mathrm{mod}\ }
\newcommand{\id}{\mathrm{id}}

\title{On real anti-bicanonical curves with one double point on the $4$-th real Hirzebruch surface. II}

\author{Sachiko Saito}

\address{Department of Mathematics Education,\ Asahikawa Campus,\\
Hokkaido University of Education,\ Asahikawa, JAPAN}
\email{saito.sachiko@a.hokkyodai.ac.jp}

\begin{document}
\pagestyle{plain}

\begin{abstract}
A real $2$-elementary K3 surfaces of type $((3,1,1),- \id)$ yields 
a real anti-bicanonical curve $s \cup  A^\prime_1$ (disjoint union) on the $4$-th real Hirzebruch surface $\bff_4$ 
where $s$ is the exceptional section of $\bff_4$ and the real curve $A^\prime_1$ has one real double point.
(See Section \ref{RealK3-311} below.) 
We give a criterion (Proposition \ref{criterion}) which determines whether the real double point is degenerate or not. 
One direction of the assertion of Proposition \ref{criterion} has already been proved in Lemma 4.6 in the preceding paper \cite{SaitoSachiko2015}. 
In this paper we prove the inverse direction. 
\end{abstract}

\maketitle

%%%%%%%%%%%%%%%%%%%%%%%%%%%%%%%%%%%%%%%%%%%%%%%%%%%%%%%%%%%%%%%%%%
\setcounter{section}{0}

\section{Introduction}

\subsection{Real $2$-elementary K3 surfaces}

In this paper we mainly discuss K3 surfaces $X$ with a non-symplectic holomorphic involution $\tau$. 
We often call them {\it $2$-elementary K3 surfaces} $(X,\tau)$ 
(\cite{Nikulin81}, \cite{AlexeevNikulin2006}, \cite{NikulinSaito05}, \cite{NikulinSaito07}, \cite{SaitoSachiko2015}, and e.t.c.). 
Note that every K3 surface with a non-symplectic holomorphic involution is algebraic. 
Hence, it has hyperplane sections. 

\begin{definition}[Real $2$-elementary K3 surface]
We say that a triple $(X,\tau,\varphi)$ is a {\it real} K3 surface with non-symplectic holomorphic involution 
(or {\it real} $2$-elementary K3 surface) 
if\\
\ \ {\rm (1)}\ $(X,\tau)$ is a K3 surface $X$ with a non-symplectic holomorphic involution $\tau$,\\
\ \ {\rm (2)}\ $\varphi$ is an anti-holomorphic involution on $X$, and\\
\ \ {\rm (3)}\ $\varphi \circ \tau = \tau \circ \varphi$.
\end{definition}

For a $2$-elementary K3 surface $(X,\tau)$, let %%%% realとは限らない
$$
{H_2}_+(X, \bz)
$$
denote the fixed part of $\tau_* : H_2(X, \bz) \to H_2(X, \bz)$. 
It is well-known that $H_2(X, \bz)$ is an even unimodular lattice of signature $(3,19)$. 
${H_2}_+(X, \bz)$ is a primitive hyperbolic $2$-elementary sublattice of $H_2(X, \bz)$. 
Note that 
$${H_2}_+(X, \bz) \subset \Pic (X),$$
where $\Pic (X)$ denotes the Picard lattice of $X$. 

Let $\bl_{K3}$ be an even unimodular lattice of signature $(3,19)$ and fix it. 
Note that the isometry class of $\bl_{K3}$ is unique. 
Let 
$$S \ \ (\subset \bl_{K3})$$
be a primitive hyperbolic $2$-elementary sublattice of $\bl_{K3}$. 

We set $r(S) := \rk S$.\ \ \ The non-negative integer $a(S)$ is defined by $S^\ast /S \cong (\bz/2\bz)^{a(S)}$. 
We define the ``parity" $\delta (S)$ of $S$ as follows. 
$$
\delta (S) := \left\{
\begin{array}{cl}
0 &\ \ \ \mbox{if}\ z \cdot \sigma (z) \equiv 0 \mmod 2 \ \ (\forall z \in \bl_{K3})\\
1 &\ \ \ \mbox{otherwise,}
\end{array}
\right.
$$

\noindent 
--------------------------------------------------------------------------------------------------------------------------------

\noindent
{\footnotesize 
partially supported by JSPS Grant-in-Aid for Challenging Exploratory Research 25610001 (2013/4 --- 2016/3).\\
{\it 2010 AMS Mathematics Subject Classification}:\ \ 14J28, 14P25, 14J10.}

\newpage

where $\sigma : \bl_{K3} \to \bl_{K3}$ is the unique integral involution whose fixed part is $S$. 

It is known that the triplet $(r(S),a(S),\delta(S))$ 
determines the isometry class of the lattice $S$ (\cite{Nikulin81}). 
Moreover, if $S$ and $S^\prime$ are isometric primitive hyperbolic $2$-elementary sublattices of the K3 lattice $\bl_{K3}$, 
then there exists an ambient automorphism 
$f$ of $\bl_{K3}$ such that 
$f(S^\prime) = S$ (\cite{AlexeevNikulin2006}, \cite{Nikulin79}). 

%%%%%%%%%%%%%%%%%%%%%%%%%%%%%%%%%%%%%%%%%%%%%%%%%%%%%%%%%%%%%%%%%%%%%%

We fix a half cone 
$$V^+(S)$$
of the cone 
$$V(S):= \{ x \in S \otimes \br\ |\ x^2 > 0\}.$$

Moreover, we fix a fundamental subdivision 
$$\Delta(S)=\Delta(S)_+\cup -\Delta(S)_+$$
of all elements with square $-2$ in $S$. 

This is equivalent to fixing a fundamental (closed) chamber (see \cite{NikulinSaito05}) 
$$\M \ \ \ \ (\subset V^+(S))$$          %%%%% \M
for the group                            %%%%% realとは限らない時点で自然なmarkingのために$\M$を指定．$\M$はeffective classesを決める．
$W^{(-2)}(S)$ 
generated by reflections in all elements with square $(-2)$ in $S$. 

Note that $\M$ and $\Delta(S)_+$ define each other by the condition $\M \cdot \Delta(S)_+ \ge 0$.
%%These additional structures $\M \subset V^+(S)$ of the hyperbolic lattice $S$ 
%%are defined uniquely up to the action of the group $\{\pm 1\}W^{(-2)}(S)$.

%%%%%%%%%%%%%%%%%%%%%% real case $\theta$ from here %%%%%%%%%%%%%%%%%%%%%%%%%%%%%%%%%%%%

Let $\theta$ be an integral involution of $S$. 

\begin{definition} \label{real_2-elementary K3_S_theta}
We say that $(X,\tau,\varphi)$ is a real $2$-elementary K3 surface {\it of type $(S,\theta)$} 
if there exists an isometry (so-called ``marking" later) 
$$\alpha : H_2(X, \bz) \cong \bl_{K3}$$
such that 
$\alpha({H_2}_+(X, \bz)) = S$ and the following diagram commutes:
$$
\begin{CD}
{H_2}_+(X, \bz) @> {\alpha}>>S\\
@V{\varphi_*}VV @VV{\theta}V\\
{H_2}_+(X, \bz) @> {\alpha}>>S .
\end{CD}
$$
\end{definition}

%%%%%%%%%%%%%%%marked real $2$-elementary K3 surface の最初の定義(by Nikulin 2005)

\begin{definition}[marked real $2$-elementary K3 surfaces] \label{marked_real_K3}
We define that a {\it marked real $2$-elementary K3 surface of type $(S,\theta)$} is 
a pair 
$$((X,\tau,\varphi),\ \alpha)$$
 of 
a real $2$-elementary K3 surface $(X,\tau,\varphi)$ of type $(S,\theta)$ (Definition \ref{real_2-elementary K3_S_theta} above) 
and an isometry, which is called {\it marking}, 
$$\alpha : H_2(X, \bz) \cong \bl_{K3}$$
such that 
\begin{itemize}
\item $\alpha({H_2}_+(X, \bz)) = S$,\\
\item $\alpha \circ \varphi_* = \theta \circ \alpha \ \ \ \text{on} \ {H_2}_+(X, \bz)$, \\
\item $\alpha_{\br}^{-1}(V^+(S))$ contains a hyperplane section of $X$, where 
$\alpha_{\br}$ stands for the real extension of $\alpha$, and \\
\item the set $\alpha^{-1}(\Delta(S)_+)$ contains only effective classes of $X$. 
\end{itemize}
\end{definition}

%% $S$と， $V^+(S)$と，subdivision $\Delta(S)=\Delta(S)_+\cup -\Delta(S)_+$と，$\theta$ をfixし，
%% 上の４条件を満たすmarkingをされているような$((X,\tau,\varphi),\ \alpha)$のみを分類する．
%% ただし，この定義のmarkingだと，$S$のどの元も指定していない．
%% あとで，(3,1,1)caseでは，特定の元を指定する．Itenbergもそうしていた．
%% 各type$S$に応じた$(X,\tau)$の幾何構造と関連させて，markingの仕方は適宜変えるわけである．
%
Note that (\cite{NikulinSaito05}) for any $(X,\tau)$, 
we can take $\alpha$ such that 
$\alpha_{\br}^{-1}(V^+(S))$ contains a hyperplane section of $X$. 

%% 必要であれば，marking $\alpha$ を $ - \alpha$に取り替える
%% つまり，最初に指定しておいた half cone V^+(S)に合わせて，K3曲面のmarkingが決まる

%%%%%%%%%%%%%次に格子の考察%%%%%%%%%%%%%%%%%%%%%%%%%%%%%%%%%%%%%%%%%%%%%%%%%%%%%%%%%%%

\subsection{Integral involutions of $\bl_{K3}$ of type $(S,\theta)$}

Let $S$ be a hyperbolic $2$-elementary sublattice of $\bl_{K3}$ and 
$\theta : S \to S$ be an integral involution (as above). 
\begin{definition}[Integral involution $\psi$ of $\bl_{K3}$ of type $(S,\theta)$]
Let $\psi : \bl_{K3} \to \bl_{K3}$ be an integral involution of the lattice $\bl_{K3}$ such that the following diagram commutes:
%%%%%%%$\psi(S)=S,\ \ {\psi |}_S = \theta$. 
$$
\begin{array}{rcl}
         S           & \subset &  \bl_{K3}          \\
\theta \ \downarrow  &         &  \downarrow \ \psi \\
         S           & \subset &  \bl_{K3} .
\end{array}
$$
We call such a pair $(\bl_{K3},\psi)$ (or $\psi$ itself) 
an {\it integral involution of $\bl_{K3}$ of type $(S,\theta)$}. 
\end{definition}

%%%%%%%%%%%%%%%%%%%%%%%%%%%%%%%%%%%%%%%%%%%%
Let $((X,\tau,\varphi),\ \alpha)$ be a marked real $2$-elementary K3 surface of type $(S,\theta)$ as above. 
If we set 
$$\psi := \alpha \circ \varphi_* \circ \alpha^{-1} : \bl_{K3} \to \bl_{K3},$$
then we have $\psi(S) = S$, and $\psi (x) = \theta (x)$ for every $x \in S$. 
Hence, $(\bl_{K3},\psi)$ is an integral involution of $\bl_{K3}$ of type $(S,\theta)$. 

We call the integral involution $(\bl_{K3},\psi)$ of type $(S,\theta)$ 
{\it the associated integral involution} with 
a marked real $2$-elementary K3 surface $((X,\tau,\varphi),\ \alpha)$ of type $(S,\theta)$ 
if the following diagram commutes:
$$
\begin{CD}
H_2(X, \bz) @>{\alpha}>>\bl_{K3} \\
@V{\varphi_*}VV @VV{\psi}V\\
H_2(X, \bz) @>{\alpha}>>\bl_{K3} .
\end{CD}
$$

%%%%%%%%%%%%%%%%%%%%%%%%%%%%%%%%%%%%
\begin{definition}[the subgroup $G$]
Let $\Delta(S,L)^{(-4)}$ be the set of all elements $\delta_1$ in $S$ such that 
$\delta_1^2=-4$ and there exists $\delta_2\in S^\perp_L$  %%%% 相棒
such that $(\delta_2)^2=-4$ and 
$\delta=(\delta_1+\delta_2)/2\in L$. 
Let $W^{(-4)}(S,L)$ 
be the subgroup of $O(S)$ generated by reflections in all elements in $\Delta(S,L)^{(-4)}$, and 
%
%%%%%$(-4)$元の直交超平面でのreflectionsはautomorphismsとなるが，
%%%%%$(-6)$元の直交超平面でのreflectionsは，ならない．
%
$W^{(-4)}(S,L)_\M$ be the stabilizer subgroup of $\M$ in $W^{(-4)}(S,L)$.
We define the subgroup 
$G$ 
to be generated by reflections $s_{\delta_1}$ in 
all elements $\delta_1\in \Delta(S,L)^{(-4)}$ which are contained either in $S_+$ or in $S_-$ and 
satisfy $(s_{\delta_1})_{\br}(\M)=\M$,   %%%%% ($\M$をstabilizeする)
where 
$s_{\delta_1}$ denotes the reflection at the orthogonal hyperplane ${\delta_1}^{\perp}$ on $S$, 
$(s_{\delta_1})_{\br}$ stands for the real extension of $s_{\delta_1}$, and 
we set 
$S_{\pm} := \{ x\in S \,|\, \theta(x) = \pm x \}$. 
Then $G$ is a subgroup of $W^{(-4)}(S,L)_\M$.
\end{definition}

\begin{definition}[Isometries with respect to the group $G$]
Let $(\bl_{K3},\psi_1)$ and $(\bl_{K3},\psi_2)$ be two integral involutions of $\bl_{K3}$ of type $(S,\theta)$. 
We define that an {\it isometry with respect to the group} $G$ from $(\bl_{K3},\psi_1)$ to $(\bl_{K3},\psi_2)$ 
is an isometry $f: \bl_{K3} \to \bl_{K3}$ 
such that $f(S)=S$, ${f|}_S \in G$, and 
the following diagram commutes:
$$
\begin{CD}
\bl_{K3} @>{f}>>\bl_{K3} \\
@VV\psi_1 V@V\psi_2 VV\\
\bl_{K3} @>{f}>>\bl_{K3} .
\end{CD}
$$
We say that 
two integral involutions $(\bl_{K3},\psi_1)$ and $(\bl_{K3},\psi_2)$ of type $(S,\theta)$ are 
{\it isometric with respect to the group $G$} 
if there exists an isometry  with respect to the group $G$ from $(\bl_{K3},\psi_1)$ to $(\bl_{K3},\psi_2)$.
By an {\it automorphism} of an integral involution $(\bl_{K3},\psi)$ of type $(S,\theta)$ {\it with respect to the group $G$} 
we mean 
an isometry with respect to the group $G$ from $(\bl_{K3},\psi)$ to itself.
Namely, 
an isometry $f: \bl_{K3} \to \bl_{K3}$ which satisfies that 
$\psi \circ f = f \circ \psi,\ f(S)=S \ \text{and} \  {f|}_S \in G$. 
\end{definition}

%%%%複素解析的同型

\begin{definition}[analytic isomorphisms with respect to $G$] \label{analytic-iso}
We say that two marked real $2$-elementary K3 surfaces 
$((X,\tau,\varphi),\ \alpha)$ and $((X^\prime,\tau^\prime,\varphi^\prime),\alpha^\prime)$ of type $(S,\theta)$ are 
{\it analytically isomorphic with respect to the group} $G$ 
if there exists an analytic isomorphism 
$f : X \to X^\prime$ such that 
$f \circ \tau = \tau^\prime \circ f$, $f \circ \varphi = \varphi^\prime \circ f$ and 
$\alpha^\prime \circ f_* \circ \alpha^{-1}|S \in G $.    %%% $f_*|S\in G$.
\end{definition}

%%%%%%%%%%%%%%%%%%%%%%%%%%%%%%%%%%%%%%%%%%%%%%%%%%%%%%%%%%%%%%%
\subsection{Period domains} \label{period domains section}

Now let us {\bf fix} an integral involution 
$$(\bl_{K3},\psi)$$
of type $(S,\theta)$ 
throughout this subsection.  %%%%% ここでは，いろいろなisometry classesを考えず，ひとつのisometry classを固定して考える．

We follow the formulations of period domains of 
marked real $2$-elementary K3 surfaces 
(see Itenberg \cite{Itenberg92} and Nikulin-Saito \cite{NikulinSaito05}). 

We set 
$$
\Omega_\psi := 
\{ \omega \ (\in \bl_{K3}\otimes \bc ) \ |\ %%%%%%% J of Singu の refereeの指摘で \neq 0 を取った．
\omega \cdot \omega =0, \ \omega \cdot \overline{\omega} >0, \ \omega \cdot S =0, \ 
\psi_{\bc}(\omega)=\overline{\omega} \}/\br^{\times} .
$$
%%%%%%%%%%%%%%%%%%%%%%%%%%%%%%%%%%%%%%%%%%%%%%%%%%%%%%%%%%%%%%%%%%%%

Let $((X,\tau,\varphi),\ \alpha)$ be a marked real $2$-elementary K3 surface of type $(S,\theta)$ 
satisfying 
$$\alpha \circ \varphi_* \circ \alpha^{-1} = \psi ,
\footnote{All marked real $2$-elementary K3 surfaces 
whose associated integral involutions are \underline{isometric} to $(\bl_{K3},\psi)$ with respect to $G$ 
satisfy $\alpha \circ \varphi_* \circ \alpha^{-1} = \psi$ if we change their markings appropriately 
(see \cite{SaitoSachiko2015}).}
$$
and let $H \ \ (\subset H_2(X,\bc))$ be the Poincare dual of $H^{2,0}(X)$. 
The $1$-dimensional subspace $\alpha_{\bc}(H)$ of $\bl_{K3}\otimes \bc$ 
is regarded as an element of $\Omega_\psi$. 

\begin{definition}
We call $\alpha_{\bc}(H)$ 
the {\it period} of a marked real $2$-elementary K3 surface $((X,\tau,\varphi),\ \alpha)$ of type $(S,\theta)$ 
satisfying $\alpha \circ \varphi_* \circ \alpha^{-1} = \psi$.
\end{definition}

\begin{definition}[Equivalence]
We say that a point $[\omega] \ (\in \Omega_\psi)$ is {\it equivalent} to a point $[\omega^\prime] \ (\in \Omega_\psi)$ 
if $[\omega^\prime] = f_{\bc}([\omega])$ for 
an automorphism $f$ of $(\bl_{K3},\psi)$ of type $(S,\theta)$ with respect to the group $G$. 
\end{definition}
%%%%%% あるK3曲面の周期だと断定できない点は，単に「point」と呼んでおく．

\begin{lemma}[\cite{SaitoSachiko2015}] \label{equivalence-lemma}
If a point $[\omega] \ (\in \Omega_\psi)$ is equivalent to $[\omega^\prime] \ (\in \Omega_\psi)$ and 
$[\omega]$ is the period of 
some marked real $2$-elementary K3 surface $((X,\tau,\varphi),\ \alpha)$ of type $(S,\theta)$ 
satisfying $\alpha \circ \varphi_* \circ \alpha^{-1} = \psi$, 
then 
$[\omega^\prime]$ is also the period of 
a marked real $2$-elementary K3 surface $((X,\tau,\varphi),\ \alpha^\prime)$ of type $(S,\theta)$ 
satisfying 
$(\alpha^\prime) \circ \varphi_* \circ (\alpha^\prime)^{-1} = \psi$ 
where $\alpha^\prime$ is another marking of $(X,\tau,\varphi)$. 
\end{lemma}

Using the global Torelli theorem, 
if two periods are equivalent, 
then corresponding marked real $2$-elementary K3 surfaces are 
analytically isomorphic (see Definition \ref{analytic-iso}). 
The converse is also true. 

The domain $\Omega_\psi$ has two connected components which are interchanged by $-\psi$. 
Since $-\psi$ is an automorphism of $(\bl_{K3},\psi)$ with respect to the group $G$, 
by Lemma \ref{equivalence-lemma}, 
it is enough to investigate the quotient space 
$$
\Omega_\psi /-\psi .
$$

Now we set 
$$\bl^\psi := \{ x \in \bl_{K3}\ |\ \psi (x) = x \},\ \ \bl_\psi := \{ x \in \bl_{K3}\ |\ \psi (x) = - x \} .$$
%%%%%% 記号変更　$\bl_{\pm} := \{ x \in \bl_{K3}\ |\ \psi (x) = \pm x \}$.

%%%%%%%%%%%%%%%%%%%%%%%%%%%% S \subset \bl_\psi case %%%%%%%%%%%%%%%%%%%%%
We restrict ourselves to the case 
$$S \subset \bl_\psi,\ \ \text{i.e.,}\ \ \theta = - \id ,$$
where ``$\id$" stands for the identity map on $S$. 
We set 
$$\bl_{-,S} := \bl_\psi \cap S^{\perp}.$$

For $[\omega] \in \Omega_\psi$ \ ($\omega \in \bl_{K3} \otimes \bc$), 
we consider the orthogonal decomposition 
$\omega = \omega_+ + i\,\omega_- \ \ (\omega_+ \in \bl^\psi \otimes \br,\ \omega_- \in \bl_\psi \otimes \br)$. 
%
% referee の指摘で \omega_+ + i \omega_- に訂正した
%
Then we have $\omega_- \in \bl_{-,S}\otimes \br$ and $\omega_+^2 = \omega_-^2 > 0$. 
Hence, both $\bl^\psi$ and $\bl_{-,S}$ are hyperbolic lattices. 
We set $V(\bl^\psi):= \{ x \in \bl^\psi \otimes \br\ |\ x^2 > 0\}$. 
$V(\bl^\psi)$ has two connected components. Let $V^+(\bl^\psi)$ be one of those (half cone). 
Let $\La_+$ denote the set of all rays (half lines) through $\zero$ in $V^+(\bl^\psi)$. 
($\La_+$ is called the {\it hyperbolic} (or {\it Lobachevsky}) space obtained from $\bl^\psi$. )
We define the hyperbolic space $\La_{-,S}$ obtained from $\bl_{-,S}$ in the same way. 

Then we have the following identification:
\begin{equation} \label{domain-quotient}
\Omega_\psi /-\psi \ = \ \La_+ \times \La_{-,S}\ \ \ \ (\mbox{a direct product}).
\end{equation}

%%%%%　{\bf BHPV} \cite{BHPV}, p.359を参照すると，$\Omega_\psi /-\psi$への全射性は言えている．
%%%% そして，ホモロジー群上のisometriesのデータから，強トレリ型定理を用いて，\\
%%%%%それに対応するようなnonsymplectic holomorphic involution$\tau$でその固定部分が$(S,\theta)$にisometricなもの，および，
%%%%%%antiholomorphic involution$\varphi$が、それぞれただ一つ存在する．\\
%%%%%つまり，対応するreal $2$-elementary K3 surface $(X,\tau,\varphi)$ of type $(S,\theta)$はある．

%%%%%%%%%%%%%%%%%%%%%%%%%%%%%%%
\subsection{($\da\br$)-nondegenerate marked real $2$-elementary K3 surfaces}

\begin{definition}
We say that an element $x (\neq 0) \in \Pic (X)\otimes \br$ is {\it nef} (for $X$) 
if 
$x \cdot C \ge 0$ for every effective curve $C$ in $\Pic (X)$. 
\end{definition}

\begin{definition}[\cite{NikulinSaito05}]\label{da-degenerate}
{\bf (i)}\ We say that 
a marked $2$-elementary K3 surface $((X,\tau),\ \alpha)$ of type $S$ is {\it ($\da$)-degenerate} 
if there exists an element $x_0 \in \alpha_{\br}^{-1}(\M)$ which is {\bf not} nef. 
Namely, $x_0 \cdot C < 0$ for an effective curve $C$ in $\Pic (X)$. 
%%%% [説明]「$\alpha^{-1}(\Delta(S)_+)$ contains only classes of effective curves of $X$」で，しかも，
%%%% all square $-2$ effective classesを含んでいるので，
%%%% 任意の$x \in \alpha_{\br}^{-1}(\M)$は，$\alpha^{-1}(\Delta(S)_+)$の任意のeffective class $C$ に対して$x \cdot C \ge 0$，
%%%%%しかも，これがいえたら，$\alpha^{-1}(S)$の任意のeffective class $C$に対して}$x \cdot C \ge 0$なのだが，
%%%%% $x_0$は，{\bf あるeffective curve $C$ in $\Pic (X)$に対して}$x_0 \cdot C < 0$となってしまう．[説明終]
%%%%%
This condition is equivalent (see \cite{Nikulin86}, \cite{AlexeevNikulin2006}) 
to 
the existence of an irreducible $(-2)$-curve on the quotient surface $Y:=X/\tau$. 
And this condition is also equivalent to 
the existence of an element $\delta \in \Pic (X)$ with $\delta^2=-2$ such that 
$\delta=(\delta_1+\delta_2)/2$
where 
$\delta_1\in \alpha^{-1}(S), \ \ \delta_2\in \alpha^{-1}(S)^\perp_{\Pic (X)}, \ \ \text{and} \ \ \delta_1^2=\delta_2^2 = -4$. 
{\bf (ii)}\ We say that 
a marked real $2$-elementary K3 surface $((X,\tau,\varphi),\ \alpha)$ of type $(S,\theta)$ is {\it ($\da\br$)-degenerate} if 
there exists a ``real" element $x_0 \in \alpha_{\br}^{-1}(S_- \cap \M)$ which is {\bf not} nef, 
where we set 
$S_{\pm} := \{ x\in S \,|\, \theta(x) = \pm x \}$.
This condition is equivalent to 
the existence of an element $\delta \in \Pic (X)$ with $\delta^2=-2$ such that 
$\delta=(\delta_1+\delta_2)/2$ 
where 
$\delta_1\in \alpha^{-1}(S),\ \  \delta_2\in \alpha^{-1}(S)^\perp_{\Pic (X)}\ \ \text{and} \ \ \delta_1^2=\delta_2^2=-4,$ 
and 
$\delta_1$ is orthogonal to an element $x \in \alpha_{\br}^{-1}(S_-\cap \text{int}(\M))$, 
where $\text{int}(\M)$ denote the interior part of $\M$, i.e., 
the polyhedron $\M$ without its faces.
%%%%%%　$\M$のface上にある点でも，$V^+(S)$においては内点である．
\end{definition}

Considering associated integral involutions, 
we have:
%%%%%%%%%%%%%%%%%%%%%%%%%%%%%%%%%%%%%%%%%%%%%%%%%%%%%%%%%%%%%%%%%%%%%
\begin{theorem}[\cite{NikulinSaito05}]\label{theorem2005moduli}
The natural map gives a bijective correspondence between 
the connected components of the period domain of 
($\da\br$)-nondegenerate marked real $2$-elementary K3 surfaces of type $(S,\theta)$ 
and 
the isometry classes with respect to $G$ of integral involutions of $\bl_{K3}$ of type $(S,\theta)$ 
such that the fixed part $\bl^\psi$ of $\psi$ is hyperbolic. 
\end{theorem}

%%%%%%%%%%%%% われわれの具体例　marked real $2$-elementary K3 surface of type $((3,1,1),- \id)$}
\section{Real $2$-elementary K3 surfaces of type $((3,1,1),- \id)$}  \label{RealK3-311}

Now we fix a sublattice $S$ of the K3 lattice $\bl_{K3}$ with the invariants
$$(r(S),a(S),\delta(S)) = (3,1,1).$$
We consider $2$-elementary K3 surfaces of type $S \ (\cong (3,1,1))$ (\cite{NikulinSaito07},\cite{SaitoSachiko2015}). 
We quote the following results from Alexeev and Nikulin \cite{AlexeevNikulin2006}. 
See also \cite{NikulinSaito07} and \cite{SaitoSachiko2015}. 

Let $(X,\tau)$ be a $2$-elementary K3 surface of type $S \cong (3,1,1)$. 
Let 
$A := X^{\tau}$
be the fixed point set (nonsingular complex curve) of $\tau$. 
Then we have 
$$A = A_0 \cup A_1\ \ \mbox{(disjoint union)},$$
where $A_0$ is a nonsingular rational curve ($\cong \bp^1$), %%%% with $A_0^2=-2$ 種数公式より自動的に
and \ $A_1$ is a nonsingular curve of genus $9$. 

$(X,\tau)$ has a structure of an elliptic pencil $|E+F|$, and 
$\tau$ is the inversion map of the group structure of the elliptic pencil 
with the zero section $A_0$. %%%% $A_0$はfiberの群構造のzeroの部分

The unique reducible fiber $E+F$ having the following properties: 
\begin{description}
\item[(i)] $E$ is a nonsingular rational curve ($\cong \bp^1$) %%%% with $E^2=-2$ 種数公式より自動的に
and $E\cdot A_0=1$.

\item[(ii)] $E\cdot F=2$, \ $F^2=-2$, \ $F\cdot A_0=0$, \ and \ 
$F$ is either \\
a nonsingular rational curve (``{\it type IIa case}"), or \\
the union of two nonsingular rational curves $F^\prime$ and $F^{\prime\prime}$ 
which are conjugate by $\tau$, $F^\prime \cdot F^{\prime\prime} = 1$ (``{\it type IIb case}"). 
%%%%%%$(F^\prime)^2 = (F^{\prime\prime})^2 = -2$

\item[(iii)] The classes $[A_0]$,\ $[E]$ and $[F]$ generate the lattice ${H_2}_+(X, \bz) \ (\cong S)$. 
Moreover, 
$A_1 \cdot E = 1,\ \ A_1 \cdot F = 2$.
The Gram matrix of the lattice ${H_2}_+(X, \bz)$ 
with respect to the basis $[E],\ [F]$ and $[A_0]$ is 
$$
\begin{array}{crrr}
        & {[E]}   & {[F]} & {[A_0]} \\
{[E]}   &   -2    &   2   &   1     \\
{[F]}   &    2    &  -2   &   0     \\
{[A_0]} &    1    &   0   &  -2 .
\end{array}
$$
\end{description}

Then we have an orthogonal decomposition
$$
{H_2}_+(X, \bz) = \bz ([A_0],\ [E]+[F]) \oplus \bz ([F]) ,
$$
where the subgroups $\bz ([A_0],\ [E]+[F])$ and $\bz ([F])$ are isometric to 
the hyperbolic plane and $\langle -2 \rangle$ respectively. 

%%%%%%%%% The quotient surface $Y$ %%%%%%%%%%%%%%%%%%%%%%%%%%%%
We now consider the quotient surface $Y := X/\tau$ (so-called ``DPN surface" (\cite{NikulinSaito05})) 
and let $\pi : X \to Y$ be the quotient map. 
We define the curves on $Y$ as follows: 
$$e:=\pi(E) \ \ \ \text{and} \ \ \ f:= \pi(F) .$$
%%%
If $F$ is the union of two nonsingular rational curves $F^\prime$ and $F^{\prime\prime}$ 
which are conjugate by $\tau$ and $F^\prime \cdot F^{\prime\prime} = 1$, 
then we have 
$f = \pi(F) = \pi(F^\prime \cup F^{\prime\prime}) = \pi(F^\prime) = \pi(F^{\prime\prime}).$
%%%% f : an irreducible curve%%%

We use the same symbols $A_0$ and $A_1$ for their images in $Y$ by $\pi$. 
Then, 
the Picard group $\Pic(Y)$ of $Y$ is generated by the classes $[e]$, $[f]$ and $[A_0]$. 
The Gram matrix of $\Pic(Y)$ with respect to the basis $[e],\ [f]$ and $[A_0]$ is 
$$
\begin{array}{crrr}
        &  {[e]}  &  {[f]}  & {[A_0]} \\
{[e]}   &    -1   &    1    &    1    \\
{[f]}   &     1   &   -1    &    0    \\
{[A_0]} &     1   &    0    &   -4 .
\end{array}
$$

We have 
$A_1\cdot e =  1 \ \ \ \mbox{and} \ \ \ A_1\cdot f = 2.$

We next contract the exceptional curve 
$f = \pi(F)$ 
to a point. Then we get a blow up
$$\BL : Y \to \bff_4,$$
where $\bff_4$ is the $4$-th Hirzebruch surface. (See \cite{NikulinSaito07}.)

We set
$s:=\BL(A_0),\ A^\prime_1 := \BL(A_1),\ c:=\BL(e)$. 
Then 
$s$ is the exceptional section of $\bff_4$ with $s^2 = -4$, and 
$c$ is a fiber of the fibration $\bff_4\to s$ with $c^2=0$. 

We have
$$\BL(A) = \BL(A_0) + \, \BL(A_1) =  s +  A^\prime_1 \ \ \in |-2K_{\bff_4}|.$$
Namely, $\BL(A)$ is \underline{an anti-bicanonical curve on $\bff_4$}. 
Since $s \cdot A^\prime_1 = 0$, 
$A^\prime_1$ does not intersect the section $s$. 
Since $-2K_{\bff_4} \sim 12c+4s$, 
we have $A^\prime_1 \ \ \in |12c+3s|$. 

%%%%%%%%%%%%ここからreal case %%%%%%%%%%%%%%%%%%%%%%%%%%%%%%%%%

\begin{remark}[\cite{NikulinSaito07}]
For any real $2$-elementary K3 surface $(X,\tau,\varphi)$ of type $(S,\theta)$ with $S \cong (3,1,1)$, 
we have 
$$\theta  =  - \id ,$$
and 
$$G = \{ \id \}.$$
It is known that all real $2$-elementary K3 surfaces of type $((3,1,1), - \id)$ 
are ($\da\br$)-non-degenerate. 
\end{remark}

\begin{remark}[\cite{NikulinSaito07}, \cite{SaitoSachiko2015}] \label{the real double point}
$F$ is a nonsingular rational curve if and only if 
$A_1$ intersects $f$ in two distinct points, and 
$F$ is a union of two nonsingular rational curves if and only if 
$A_1$ touches $f$. 

If $A_1$ intersects with $f$ at two distinct points, then they are real points or non-real conjugate points. 
In the former case the real double point of $A^\prime_1$ is a real node, and 
in the latter case it is a real isolated point. 
Anyway the real double point of $A^\prime_1$ is nondegenerate. 
If $A_1$ touches to $f$ in $Y$, then the real double point is a real cusp (a degenerate double point). 
\end{remark}

Since all real $2$-elementary K3 surfaces of type $((3,1,1), - \id)$ are ($\da\br$)-non-degenerate, 
by Theorem \ref{theorem2005moduli}, 
the connected components of the period domain of 
marked real $2$-elementary K3 surfaces of type $((3,1,1), - \id)$ 
and 
the isometry classes with respect to $G = \{ \id \}$ 
of integral involutions of $\bl_{K3}$ of type $((3,1,1), - \id)$ such that 
the fixed part $\bl^\psi$ of $\psi$ is hyperbolic 
are in bijective correspondence. 

However, as is written in Remark 8 in Subsection 2.4 in \cite{NikulinSaito07}, 
it is possible that 
an isometry class of integral involutions of $\bl_{K3}$ of type $((3,1,1), - \id)$ 
corresponds to {\bf both type IIa case and type IIb case}. 
In other words, such isometry classes with respect to $G$ (in the sense of Theorem \ref{theorem2005moduli}) cannot distinguish 
{\bf degenerate and nondegenerate double points} of the curves $A^\prime_1$ on $\bff_4$ 
(Recall Remark \ref{the real double point}). 

Therefore, we would like to determine whether the real double point of the curve $A^\prime_1$ on $\bff_4$ is degenerate or not. 
See also Lemma 4.6 and Problem 1 in \cite{SaitoSachiko2015}. 
In order to do this, 
we define more strict markings of real $2$-elementary K3 surfaces of type $((3,1,1),- \id)$. 
Compare the following Definition \ref{marked311} to the Definition \ref{marked_real_K3} above. 

%%%%%%%記号の準備
Let 
$
\Ea,\ \ \Fa,\ \ \mbox{and}\ \ \Aa
$
be generators of $S$ with the Gram matrix 
$$
\begin{array}{crrr}
       & {\Ea}  & {\Fa} & {\Aa}  \\
{\Ea}  &   -2   &   2   &     1 \ \\
{\Fa}  &    2   &  -2   &     0 \ \\
{\Aa}  &    1   &   0   &    -2
\end{array}
$$
where 
$\M \cdot \Ea \ge 0$,\ 
$\M \cdot \Fa \ge 0$,\ and \ 
$\M \cdot \Aa \ge 0$.   %%% (effectiveになるように)

And we set 
$$\bu := \bz (\Aa,\ \Ea + \Fa).$$
$\bu$ is isometric to the hyperbolic plane, and 
$$
S = \bu \oplus \bz (\Fa)
$$
is an orthogonal decomposition of $S$. 

%% $(3,1,1)$の場合に対して，もっと厳格な marked real $2$-elementary K3 surface の定義\\
%% 以前の定義 \ref{marked_real_K3}より精密である．結局，生成元のmarkingによる行先をすべて指定．\\
%% Itenberg の$<2> \oplus <-2>$の場合のmarkingも，元まで指定している．\\
%% 各$(S,\theta)$によって幾何的状況($A_1$ と $f$ の位置関係など)を特徴づける元が違うため，それらを記述するために，
%% 「特定の元」に関する条件をつける．\\

\begin{definition}[marked real $2$-elementary K3 surfaces $((X,\tau,\varphi), \alpha)$ of type $((3,1,1),- \id)$] \label{marked311}
We define that a {\it marked real $2$-elementary K3 surfaces $((X,\tau,\varphi), \alpha)$ of type $((3,1,1),- \id)$} 
is 
a pair of a real $2$-elementary K3 surface $(X,\tau,\varphi)$ of type $((3,1,1),- \id)$ (Definition \ref{real_2-elementary K3_S_theta}) 
and a marking (isometry) 
$$\alpha : H_2(X, \bz) \cong \bl_{K3}$$
such that 
\begin{itemize}
\item $\alpha({H_2}_+(X, \bz)) = S$,
\item $\alpha \circ \varphi_* = \theta \circ \alpha \ \text{on} \ {H_2}_+(X, \bz)$,
\item $\alpha_{\br}^{-1}(V^+(S))$ contains a hyperplane section of $X$,
\item the set $\alpha^{-1}(\Delta(S)_+)$ contains only classes of effective curves of $X$, 
%%%%%%%%%%%%%%%%%%%%%%%%%%%%%%%
\item $\alpha([A_0]) = \Aa$,\ $\alpha([E]) = \Ea$, and $\alpha([F])=\Fa$.  %%%%各元の行先を指定
\end{itemize}
\end{definition}

Note that every real $2$-elementary K3 surface $(X,\tau,\varphi)$ of type $((3,1,1),- \id)$ 
has such a marking $\alpha$. 

We give a criterion for 
the unique double point of a real anti-bicanonical curve $\BL(A) = s +  A^\prime_1 $ on $\bff_4$ 
with one real double point on $A^\prime_1$ 
to be nondegenerate. 

Let us consider and fix an integral involution 
$$(\bl_{K3},\psi)$$
of type $((3,1,1),- \id)$ again, %%%%% ひとつ固定
and consider the period domain \eqref{domain-quotient} in Subsection \ref{period domains section} 
$$
\Omega_\psi /-\psi \ = \ \La_+ \times \La_{-,S} \ .
$$

\begin{proposition} \label{criterion}
Let $[[\omega]]$ be a point
\footnote{
%%%%%%% Remark \ref{remove-hyperplanes}の直交超平面たちを除かなくても全射性は言えている．
%%%%%%% しかもその対応するK3曲面が求める性質を満たす．
By the surjectivity of the period map of marked K3 surfaces (\cite{BHPV}, p.339), 
$[\omega]$ is the period of a marked real $2$-elementary K3 surface $((X,\tau,\varphi), \alpha)$ of type $((3,1,1),- \id)$ 
satisfying $\alpha \circ \varphi_* \circ \alpha^{-1} = \psi$.}
 in $\Omega_\psi /-\psi$. 
Then, the real double point of the curve $A^\prime_1$ on $\bff_4$ (Remark \ref{the real double point}) 
which corresponds to $((X,\tau,\varphi), \alpha)$ 
is nondegenerate 
if and only if 
there are no $\bv \ (\neq \pm \Fa)$ in $\bl_{K3}$ satisfying:
$$
\bv \cdot \omega = 0,\ \ \bv \cdot \bu =0,\ \ \text{and}\ \ \bv^2 = -2.
$$
\end{proposition}

\begin{proof}
The ($\Longleftarrow$) direction has been proved in Lemma 4.6 of \cite{SaitoSachiko2015}. 
Here we prove the ($\Longrightarrow$) direction. 
Assume that the double point of the real anti-bicanonical curve $\BL(A)$ on $\bff_4$ is nondegenerate. 
Then $F$ is irreducible (\cite{NikulinSaito07}, \cite{SaitoSachiko2015}). 
Suppose that there {\bf exists} a $\bv \ (\neq \pm \Fa)$ in $\bl_{K3}$ satisfying that 
$\bv \cdot \omega = 0,\ \bv \cdot \bu =0$, and $\bv^2 = -2$. 
Then $-\bv$ has the same properties. 
By Riemann-Roch Theorem, we have $l(\alpha^{-1}(\bv))+l(-\alpha^{-1}(\bv)) \geq \alpha^{-1}(\bv)^2/2 +2 = 1$. 
Hence, $\alpha^{-1}(\bv)$ or $-\alpha^{-1}(\bv)$ is effective. (See \cite{AlexeevNikulin2006}, p.23 or \cite{BHPV}, p.311.) 
Hence, we may assume $\alpha^{-1}(\bv)$ is an effective class. 
Then, $\tau^*(\alpha^{-1}(\bv))$ is also effective. 
We have $\tau^*(\alpha^{-1}(\bv))^2 = -2$, $\tau^*(\alpha^{-1}(\bv))\cdot [A_0] = 0,\ \tau^*(\alpha^{-1}(\bv))\cdot ([E]+[F])=0$ 
and $\alpha(\tau^*(\alpha^{-1}(\bv)))\cdot \omega = 0$ since $\tau^*$ is non-symplectic. 
Hence, $\alpha(\tau^*(\alpha^{-1}(\bv)))$ has the same properties as $\bv$. 
Since $\alpha(\tau^*(\alpha^{-1}(\bv))) + \bv \in S$, it is contained in $\bz (\Fa)$. 
There exists an integer $n$ such that 
$$\alpha(\tau^*(\alpha^{-1}(\bv))) + \bv = n\Fa.$$
Let us consider a hyperplane section class $h \ \in \Pic X$, which is ample. 
Since $\tau^*(\alpha^{-1}(\bv))$ and $\alpha^{-1}(\bv)$ are effective, 
we have $\tau^*(\alpha^{-1}(\bv))\cdot h >0$ and $\alpha^{-1}(\bv)\cdot h >0$. (Nakai's criterion) 
We have $\tau^*(\alpha^{-1}(\bv)) + \alpha^{-1}(\bv) \neq 0$. 
%%%Actually, if $\tau^*(\alpha^{-1}(\bv)) + \alpha^{-1}(\bv) =0$, then 
%%%$(\tau^*(\alpha^{-1}(\bv)) + \alpha^{-1}(\bv))\cdot h =0$. This is a contradiction. 
Hence, we have $n \neq 0$. 
Moreover, since $\tau^*(\alpha^{-1}(\bv)) + \alpha^{-1}(\bv)$ and $[F] (= \alpha^{-1}(\Fa))$ are effective, 
we have \underline{$n>0$}. 
Actually, if $n<0$, then $0=(\tau^*(\alpha^{-1}(\bv)) + \alpha^{-1}(\bv) + (-n)\alpha^{-1}(\Fa))\cdot h >0$. This is a contradiction. 

If $\tau^*(\alpha^{-1}(\bv)) = \alpha^{-1}(\bv)$, then $(\tau^*(\alpha^{-1}(\bv)) + \alpha^{-1}(\bv))^2 = (2\alpha^{-1}(\bv))^2 = -8$. 
On the other hand, $(\tau^*(\alpha^{-1}(\bv)) + \alpha^{-1}(\bv))^2 = (n[F])^2 = -2n^2$. 
Hence, $n=2$. Thus we have $\bv = \Fa$. This is a contradiction. 
Thus we have $\tau^*(\alpha^{-1}(\bv)) \neq \alpha^{-1}(\bv)$. 

Let $\alpha^{-1}(\bv)=\displaystyle \sum_i m_i \bv_i$, where $m_i$ are positive integers and $\bv_i$ are irreducible effective classes. 
Obviously we have $\alpha(\bv_i) \cdot \omega = 0$ for any $i$. %%%Picに入るから

If $\bv_i = [E]$ or $[F]$, then $\bv_i \cdot ([E]+[F]) = 0$. 
And if $\bv_i \neq [E]$ and $\neq [F]$, 
then we have $\bv_i \cdot ([E]+[F]) \geq 0$ for any $i$. %%%% irreducible同士でも，交わらなければ，交点数0
Here $E$ is irreducible, and moreover, $F$ is irreducible by the assumption. 
Since $\alpha^{-1}(\bv)\cdot ([E]+[F])=0$, we have $\bv_i \cdot ([E]+[F]) = 0$ for any $i$. 
Since $[A_0] \cdot ([E]+[F]) = 1$, we see $\bv_i \neq [A_0]$ for any $i$. 
Thus we also have $\bv_i \cdot [A_0] \geq 0$. 
Since $\alpha^{-1}(\bv)\cdot [A_0]=0$, we have $\bv_i \cdot [A_0] = 0$ for any $i$. 

Thus, we have $\alpha(\bv_i) \cdot \bu = 0$ for any $i$. 
Since $\bu$ is of signature $(1,1)$, we have $(\bv_i)^2 <0$ for any $i$ by the Hodge index theorem. 
%%Now we have 
%%$$-2 = (\alpha^{-1}(\bv))^2 = (\sum_i m_i \bv_i)^2 = \sum_i (m_i)^2(\bv_i)^2 + 2\sum_{i<j} m_i m_j \bv_i \bv_j .$$
%%Hence, there {\bf exist} $\bv_i$'s with negative squares. 
For such $\bv_i$'s, we have $\bv_i^2 = -2$,           %%% see [AlexNiku] p.22の種数公式から
and $\bv_i$'s are nonsingular rational ($\cong \bp^1$). 

Suppose $\alpha(\bv_i) = \Fa$ for some $i$, say $i=1$. Namely, $\bv_1 = [F]$. 
Let $\bv^\prime := \alpha^{-1}(\bv) - m_1 [F]$. %%%% \Fa 成分をとった残り
Since $\tau^*(\alpha^{-1}(\bv)) \neq \alpha^{-1}(\bv)$ (see above), we have 
$\bv^\prime \neq 0$.     %%%もしそうでないなら，= になる．
Thus, there exists a $\bv_i (\neq \pm [F])$ in $H_2(X,\bz)$ such that 
$\bv_i \cong \bp^1$ (irreducible), $\bv_i^2 = -2$, $\alpha(\bv_i) \cdot \omega = 0$, and $\alpha(\bv_i) \cdot \bu = 0$. 
%% irreducibleにとれたことが重要！

%%Then we have 
%%$$\alpha(\tau^*(\alpha^{-1}(\bv^\prime))) + \bv^\prime = (n - 2m_1)\Fa .$$
%%$Since $\alpha(\tau^*(\alpha^{-1}(\bv^\prime)))$ and $\bv^\prime$ are effective and nonzero, we have $n - 2m_1$ is a positive integer. 
%%% ampleなhとのintersection numberを考えるとわかる

Eventually we can choose $\bv$ such that $\alpha^{-1}(\bv)$ is an irreducible class. 

Since $\tau^*(\alpha^{-1}(\bv)) \neq \alpha^{-1}(\bv)$, 
$\tau^*(\alpha^{-1}(\bv))$ and $\alpha^{-1}(\bv)$ are represented by different irreducible curves respectively. 
Hence, we have $\tau^*(\alpha^{-1}(\bv)) \cdot \alpha^{-1}(\bv) \geq 0$. 

Since 
$(\alpha(\tau^*(\alpha^{-1}(\bv))) + \bv)^2 = (-2) + (-2) + 2(\tau^*(\alpha^{-1}(\bv)) \cdot \alpha^{-1}(\bv)) = -2n^2$, 
we have $2-n^2 \geq 0$. 
Thus, we have 
$n=1$. Namely, we have 
$$\Fa = \bv + \alpha(\tau^*(\alpha^{-1}(\bv))),$$
i.e.,
$[F] = \alpha^{-1}(\bv) + \tau^*(\alpha^{-1}(\bv))$.
Let $F^\prime$ be an irreducible curve representing $\alpha^{-1}(\bv)$, and we set 
$F^{\prime \prime} := \tau^*(F^\prime)$. 
Thus $F^\prime + F^{\prime \prime}$ represents the class $[F]$. 
Hence, there exists a marked real $2$-elementary K3 surface corresponding to the period $[\omega]$ 
which has 
$E + F^\prime + F^{\prime \prime}$ 
as the unique reducible fiber of its elliptic fibration. %%%%($\widetilde{\mathbb{A}_2}$ type)
%%%%% (そのようなElliptic fibrationからK3を作ればよい．)
Conversely, every marked real $2$-elementary K3 surface corresponding to the period $[\omega]$ 
has $E + F^\prime + F^{\prime \prime}$ as the unique reducible fiber of its elliptic fibration. 
%%%
Hence, we have $F= F^\prime + F^{\prime \prime}$. 
This contradicts the assumption that $F$ is irreducible. 
This completes the proof of Proposition \ref{criterion}.
\end{proof}

\begin{remark}
For curves of degree $6$ with one double point on $\br P^2$, 
the corresponding criterion to Proposition \ref{criterion} is written 
on the top of p. 281 in Itenberg's paper \cite{Itenberg92}. 
\end{remark}

As is written in \cite{SaitoSachiko2015}, 
we now get 
the precise image ($\subset \La_+ \times \La_{-,S}$) of the period map 
on the set of all marked real $2$-elementary K3 surfaces $((X,\tau,\varphi), \alpha)$ of type $((3,1,1),- \id)$ 
for which the real double points of the curves $A^\prime_1$ on $\bff_4$ are {\bf nondegenerate}. 
Thus we are able to continue the interesting arguments as in \cite{Itenberg92}. 

\noindent
{\bf Acknowledgments.} The author would like to thank 
Professor Hisanori Ohashi for his helpful suggestions and comments, 
Professor Ilia Itenberg for his inspiring papers, and 
Professor Viacheslav Nikulin for his constant encouragement. 

%%%%%%%%%%%%% 残りと続きは，ファイルlemma_proof2015Sep-New-tsuzuki.tex に移動させた．

%%%%%%%%%%%%%%%%%%%%%%%%%%%%%%%%%%%%%%%%%%%%%%%%%%%%%%

\end{document}